\documentclass[11pt]{article}
\usepackage{amscd, amsmath, amssymb, graphicx, floatflt}
\usepackage[all,cmtip]{xy}

\title{Line bundles with partially vanishing cohomology}
\author{Burt Totaro}
\date{  }

\def\Z{\text{\bf Z}}
\def\Q{\text{\bf Q}}
\def\R{\text{\bf R}}

\def\P{\text{\bf P}}
\def\F{\text{\bf F}}
\def\arrow{\rightarrow}

\def\imp{\Rightarrow}

\def\qed{\ QED }

\def\reg{\text{reg}}
\def\dim{\text{dim}}

\def\Hom{\text{Hom}}
\def\Tor{\text{Tor}}
\def\Spec{\text{Spec}}
\def\op{\text{op}}
\def\RR{{\mathcal R}}
\def\X{\widetilde{X}}
\def\Tot{\text{Tot}}
\def\Ext{\text{Ext}}
\def\HHom{{\mathcal Hom}}

\begin{document}
\maketitle

\newtheorem{theorem}{Theorem}[section]
\newtheorem{corollary}[theorem]{Corollary}
\newtheorem{lemma}[theorem]{Lemma}
\newtheorem{conjecture}[theorem]{Conjecture}
\newtheorem{question}[theorem]{Question}

Ample line bundles are fundamental to algebraic geometry.
The same notion of ampleness arises in many ways: geometric
(some positive multiple gives a projective embedding),
numerical (Nakai-Moishezon, Kleiman), or cohomological
(Serre) \cite[Chapter 1]{Lazarsfeld}.
Over the complex numbers, ampleness of a line bundle is also equivalent
to the existence of a metric with positive curvature (Kodaira).

The goal of this paper is to study weaker notions of ampleness,
and to prove some of the corresponding equivalences.
The subject began with Andreotti-Grauert's theorem 
that on a compact complex manifold $X$ of dimension $n$, a hermitian
line bundle $L$ whose curvature form has at least $n-q$ positive
eigenvalues at every point has $H^i(X,E\otimes L^{\otimes m})=0$
for every $i>q$, every coherent sheaf $E$ on $X$, and every integer
$m$ at least equal to some $m_0$
depending on $E$ \cite{AG}. Call the latter
property {\it naive $q$-ampleness }of $L$, for a given natural number
$q$. Thus naive $0$-ampleness is the usual notion
of ampleness, while every line bundle is naively $n$-ample. Nothing is known
about the converse to Andreotti-Grauert's theorem for $q>0$, but we can still
try to understand naive $q$-ampleness for projective varieties over any field.

Sommese gave a clear
geometric characterization of naive $q$-ampleness when in addition $L$
is semi-ample (that is,
some positive multiple of $L$ is spanned by its global sections).
In that case, naive $q$-ampleness is equivalent
to the condition that the morphism to projective space
given by some multiple of $L$ has fibers of dimension at most $q$
\cite{Sommese}. That has been useful, but
the condition of semi-ampleness is restrictive, and in this paper
we do not want to assume it.
For example, the line bundle $O(a,b)$ on $\P^1\times \P^1$ is naively 1-ample
exactly when at least one of $a$ and $b$ is positive, whereas semi-ampleness
would require both $a$ and $b$ to be nonnegative. Intuitively,
$q$-ampleness means that a line bundle is positive
``in all but at most $q$ directions''. One can hope to relate $q$-ampleness
to the geometry of subvarieties of intermediate dimension.
\begin{figure*}
\centering
\includegraphics[width=0.7\textwidth]{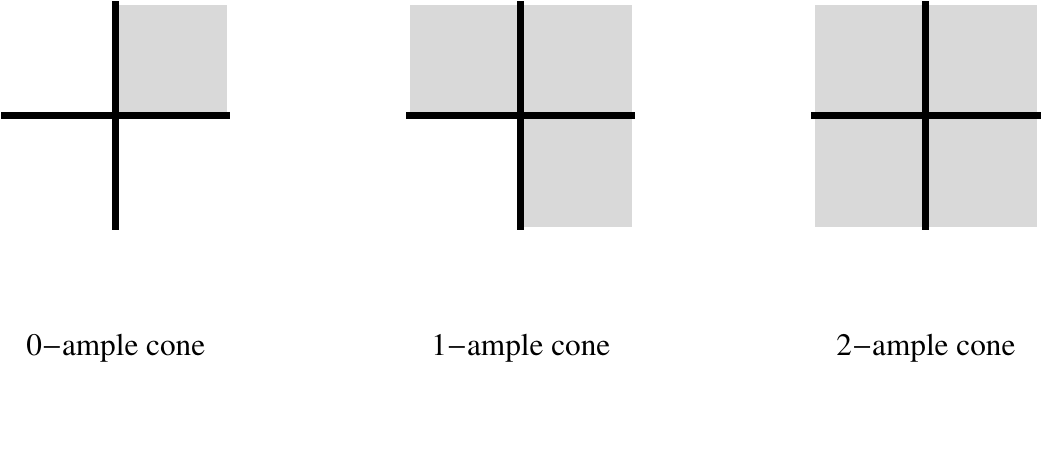}
\caption{The (open) naive $q$-ample cones in $N^1(\P^1\times \P^1)\cong \R^2$}
\end{figure*}

The main results of this paper apply to projective
varieties over a field of characteristic
zero. In that situation, we show that naive $q$-ampleness (which
is defined using the vanishing of infinitely many cohomology groups)
is equivalent to the vanishing of finitely many cohomology groups,
a condition we call $q$-T-ampleness (Theorem \ref{equivred}).
A similar equivalence holds for all projective schemes over
a field of characteristic zero (Theorem \ref{equiv}). These
results are analogues
of Serre's characterization of ampleness. Indeed, $q$-T-ampleness
is an analogue of the geometric
definition of ampleness by some
power of $L$ giving a projective embedding; the latter is also
a ``finite'' condition, unlike the definition of naive $q$-ampleness.
The equivalence implies in particular that
naive $q$-ampleness is Zariski open on families of varieties
and line bundles in characteristic zero, which is not at all
clear from the definition.

Theorem \ref{equivred} also shows that naive $q$-ampleness
in characteristic zero is equivalent to uniform $q$-ampleness,
a variant defined by Demailly-Peternell-Schneider \cite{DPS}.
It follows that naive $q$-ampleness defines an open cone (not necessarily
convex) in the N\'eron-Severi vector space
$N^1(X)$. (For example, the $(n-1)$-ample cone of an $n$-dimensional
projective variety is the complement
of the negative of the closed effective cone, by Theorem \ref{n-1}.)
After these results, it makes sense to say simply
``$q$-ample'' to mean any of these equivalent notions for
line bundles in characteristic zero.

The following tools are used for these equivalences.
First, we use the relation found by Kawamata between
Koszul algebras and resolutions of the diagonal (Theorem \ref{orlov}).
Following Arapura \cite[Corollary 1.12]{ArapuraPart},
we show that Castelnuovo-Mumford
regularity behaves well under tensor products on any
projective variety (Theorem \ref{product}). We prove
the vanishing
of certain Tor groups associated to the Frobenius homomorphism
on any commutative $\F_p$-algebra (Theorem \ref{flat}). Finally,
Theorem \ref{van} generalizes
Arapura's positive characteristic vanishing theorem
\cite[Theorem 5.4]{ArapuraPart} to singular varieties.
The main equivalences of the paper, which hold in characteristic
zero, are proved by the unusual method of reducing modulo
many different prime numbers and combining the results.

Finally, we give a counterexample to a Kleiman-type characterization
of $q$-ample line bundles.
Namely, we define an $\R$-divisor $D$ to be {\it $q$-nef }if
$-D$ is not big on any $(q+1)$-dimensional
subvariety of $X$. The $q$-nef cone is closed in $N^1(X)$ (not necessarily
convex), and all $q$-ample line bundles
are in the interior of the $q$-nef cone. The converse would
be a generalization of Kleiman's numerical criterion for ampleness
(the case $q=0$). This ``Kleiman criterion'' is true for $q=0$
and $q=n-1$, but section \ref{q-nef} shows
that it can fail for 1-ample line bundles
on a complex projective 3-fold.

Thanks to Anders Buch and Lawrence Ein for useful discussions.

\section{Notation}
\label{notation}
Define an $\R$-divisor on a projective variety $X$ to be
an $\R$-linear combination of Cartier divisors on $X$ \cite{Lazarsfeld}.
Two $\R$-divisors
are called numerically equivalent if they have the same intersection
number with all curves on $X$. We write $N^1(X)$ for the real vector space
of $\R$-divisors modulo numerical equivalence, which has finite dimension.

The {\it closed effective cone }is the closed convex cone in $N^1(X)$ spanned
by effective Cartier divisors. An line bundle is {\it pseudoeffective }if
its class in $N^1(X)$ is in the closed effective cone. A line bundle is
{\it big }if its class in $N^1(X)$ is in the interior of the closed
effective cone. A line bundle $L$ is big if and only if there are constants
$m_0$ and $c>0$ such that $h^0(X,L^{\otimes m})\geq cm^n$
for all $m\geq m_0$, where $n$ is the dimension of $X$
\cite[v.~1, Theorem 2.2.26]{Lazarsfeld}.

Let $X$ be a projective scheme of dimension $n$ over a field. Assume that
the ring $O(X)$ of regular functions on $X$ is a field; this holds
in particular if $X$ is connected and reduced. Write $k=O(X)$.
Define a line bundle $O_X(1)$ on $X$
to be {\it $N$-Koszul, }for a natural number $N$,
if the sections of $O_X(1)$ give a projective
embedding of $X$ and the homogeneous
coordinate ring $A=\oplus_{j\geq 0} H^0(X,O(j))$ is $N$-Koszul.
That is, the field $k$ has a resolution as a graded $A$-module,
$$\cdots \arrow M_1\arrow M_0\arrow k\arrow 0,$$
with $M_i$ a free module
generated in degree $i$ for $i\leq N$. (Note that
Polishchuk-Positselski's book on Koszul algebras
uses ``$N$-Koszul'' in a different sense \cite[section 2.4]{PP}.)
In particular,
we say that a line bundle $O_X(1)$ is {\it Koszul-ample }if it
is $2n$-Koszul. For example, the standard
line bundle $O(1)$ on $\P^n$ is Koszul-ample.
Backelin showed that a sufficiently
large multiple of every ample line bundle on a projective variety
is Koszul-ample (actually with $M_i$ generated in degree $i$ for
all $i$, not just $i$ at most $2n$,
but we don't need that) \cite{Backelin}. Explicit bounds on
what multiple is needed have been given \cite{Pareschi,ERT,Hering}.

One advantage of working with $N$-Koszulity for finite $N$ rather
than Koszulity in all degrees is that $N$-Koszulity is a Zariski open
condition in families, as we will use in some arguments.

Given a Koszul-ample line bundle $O_X(1)$ on a projective scheme $X$,
define the Castelnuovo-Mumford {\it regularity }of a coherent
sheaf $E$ on $X$ to be the least integer $m$ such that
$$H^i(X,E(m-i))=0$$
for all $i>0$ \cite[v.~1, Definition 1.8.4]{Lazarsfeld}. (Thus the regularity
is $-\infty$ if $E$ has zero-dimensional support.)
We know
that for any coherent sheaf $E$, $E(m)$ is globally generated and
has vanishing higher cohomology for $m$ sufficiently large, and
one purpose of regularity is to estimate how large $m$ has to be.
Namely, if $\reg(E)\leq m$, then $E(m)$ is globally generated
and has no higher cohomology \cite[v.~1, Theorem 1.8.3]{Lazarsfeld}.
We remark that if ``$H^i(X,E(m-i))=0$ for all $i>0$''
holds for one value
of $m$, then it also holds for any higher value of $m$
\cite[v.~1, Theorem 1.8.5]{Lazarsfeld}. (We generalize
this in Lemma \ref{shift}.) It is immediate
from the definition that
$\reg(E(j))=\reg(E)-j$ for every integer $j$. 

To avoid confusion, note that $\reg(O_X)$ can be greater than 0,
in contrast to what happens in the classical case of $X=\P^n$ with the standard
line bundle $O(1)$.

\section{Resolution of the diagonal}
\label{koszul}

In this section, we show that Koszul-ampleness
of a line bundle $O_X(1)$ leads to an explicit
resolution of the diagonal as a sheaf on $X\times X$. The resolution
was constructed for $O_X(1)$ sufficiently ample
by Orlov \cite[Proposition A.1]{Orlov}, and under the more
convenient assumption that
the coordinate ring of $(X,O_X(1))$ is a Koszul algebra
by Kawamata \cite[proof of Theorem 3.2]{Kawamata}. Theorem \ref{orlov}
works out the analogous statement when $O_X(1)$ is only $N$-Koszul.

Let $X$ be a projective scheme over a field. Assume that
the ring $O(X)$ of regular functions on $X$ is a field; this holds
in particular if $X$ is connected and reduced. Write $k=O(X)$.
Let $O_X(1)$ be an $N$-Koszul line bundle, for a positive integer $N$.
That is, $O_X(1)$
is very ample and the homogeneous
coordinate ring $A=\oplus_{i\geq 0}H^0(X,O(i))$ is $N$-Koszul.
(In later sections, we will work with a Koszul-ample
line bundle $O_X(1)$, which means taking $N=2\,\dim(X)$.)

Define vector spaces $B_m$ inductively by
$B_0=k$, $B_1=H^0(X,O(1))$, and
\begin{equation*}
B_m=\ker(B_{m-1}\otimes H^0(X,O(1))\arrow B_{m-2}\otimes H^0(X,O(2))).
\end{equation*}
(Products are over $k$ unless otherwise specified.)
By definition of $N$-Koszulity, the complex
\begin{equation*}
B_{N}\otimes A(-N)\arrow \cdots \arrow B_1\otimes A(-1)\arrow A
\arrow k\arrow 0
\tag{1}
\end{equation*}
of graded $A$-modules is exact. (For an integer $j$ and a graded module $M$,
$M(j)$ means $M$ with degrees lowered by $j$.) The vector space $\Tor_i^A(M,k)$
for a bounded below $A$-module $M$ can be viewed as the generators
of the $i$th step of the minimal resolution of $M$. The Koszul resolution
(1) of $k$ as an $A$-module is clearly minimal, and so
$$B_m\cong \Tor_m^A(k,k)$$
for $0\leq m\leq N$.

Let $\RR_0=O_X$, and
$$\RR_m=\ker(B_m\otimes O_X\arrow B_{m-1}\otimes O_X(1))$$
for $m>0$. The definition of $B_m$ gives a complex of sheaves
\begin{equation*}
0\arrow \RR_m\otimes_{O_X}O_X(-m+1)\arrow B_m\otimes_k O_X(-m+1)
\arrow \cdots\arrow B_1\otimes_k O_X\arrow O_X(1)\arrow 0.
\tag{2}
\end{equation*}
This is exact for $0<m\leq N$ by (1), using that (by Serre)
a sequence of sheaves is exact if tensoring with $O(j)$ for $j$ large
and taking global sections gives an exact sequence.

\begin{theorem}
\label{orlov}
There is an exact sequence of sheaves on $X\times_k X$,
\begin{equation*}
\RR_{N-1}\boxtimes O_X(-N+1)\arrow\cdots\arrow \RR_1\boxtimes O_X(-1)
\arrow \RR_0\boxtimes O_X\arrow O_{\Delta}\arrow 0,
\tag{3}
\end{equation*}
where $\Delta\subset X\times_k X$ is the diagonal.
\end{theorem}

Here $E\boxtimes F$ denotes the external tensor product
$\pi_1^*(E)\otimes \pi_2^*(F)$ on a product scheme $X_1\times_k X_2$,
for sheaves $E$ on $X_1$ and $F$ on $X_2$.

{\bf Proof. }By Serre, this sequence of sheaves is exact if tensoring with
$O(j,l)$ for all $j$ and $l$ sufficiently large and taking global sections
gives an exact sequence. The definition of $\RR_m$ implies
that 
\begin{align*}
H^0(X,\RR_m(j)) &=\ker(B_m\otimes H^0(X,O(j))\arrow
B_{m-1}\otimes H^0(X,O(j+1)))\\
 &=\ker(B_m\otimes A_j\arrow B_{m-1}\otimes A_{j+1})
\end{align*}
for all $j\geq 0$. Thus we want to prove exactness of the 
complex of $k$-vector spaces
\begin{equation*}
A_{l-N+1}\otimes\ker(B_{N-1}\otimes A_j\arrow B_{N-2}\otimes A_{j+1})
\arrow \cdots \arrow A_{l-1}\otimes\ker(B_1\otimes A_j\arrow A_{j+1})
\arrow A_l\otimes A_j\arrow A_{j+l}\arrow 0\\
\tag{4}
\end{equation*}
for all $j$ and $l$ sufficiently large. In fact, we will prove
this for all $j,l\geq 0$.

Because $A$ is an associative algebra with augmentation $A\arrow k$,
the groups 
$\Ext^*_A(k,k)$ form an associative algebra (typically not graded-commutative,
even when $A$ is commutative). The product can be 
viewed as composition in the derived category of $A$. Since
$\Ext^i_A(k,k)\cong \Tor_i^A(k,k)^*$, we can also say that
$\Tor_*^A(k,k)$
is a coassociative coalgebra \cite[section 1.1]{PP}. Thus we have natural
maps $B_{i+j}\arrow B_i\otimes B_j$ for $i$ and $j$ at most $N$.
Coassociativity says in particular that 
the two compositions $B_i\arrow A_1\otimes B_{i-1}\arrow A_1\otimes
B_{i-2}\otimes A_1$ and $B_i\arrow B_{i-1}\otimes A_1\arrow A_1
\otimes B_{i-2}\otimes A_1$ are equal, where we have identified
$B_1$ with $A_1$. Also, there is a natural isomorphism $\Ext^*_A(k,k)\cong
\Ext^*_{A^{op}}(k,k)$ that reverses the order of multiplication
\cite[section 1.1]{PP}. For a graded associative algebra $A$,
it follows that $A$ is $N$-Koszul if and only if $A^{\op}$
is $N$-Koszul. Therefore, for an $N$-Koszul algebra $A$, the following
version of the Koszul complex (1) (using the maps
$B_i\arrow A_1\otimes B_{i-1}$ rather than $B_i\arrow B_{i-1}\otimes A_1$)
is also exact:
\begin{equation*}
A(-N)\otimes B_{N}\arrow \cdots \arrow A(-1)\otimes B_1\arrow A
\arrow k\arrow 0
\tag{1'}
\end{equation*}

Let us artificially define $B_i$ to be zero for $i>N$.
Consider the following triangular diagram,
where the rows are obtained from the Koszul
complex (1') and the columns from the Koszul complex (1):

\xymatrix{
0 \ar[r] &A_0\otimes B_0 \otimes A_{j+l} \ar[r] &0 &&&\\
& \vdots\ar[u] && 0 & & \\
0\ar[r] & A_0\otimes B_{l-1}\otimes A_{j+1} \ar[r]\ar[u] &
\cdots \ar[r] & A_{l-1}\otimes B_0\otimes A_{j+1} \ar[r]\ar[u] &
0 & \\
0 \ar[r]& A_0\otimes B_l\otimes A_j \ar[r]\ar[u] & \cdots 
\ar[r] &
A_{l-1}\otimes B_1\otimes A_j\ar[r]\ar[u] &
A_l\otimes A_j \ar[r]\ar[u]& 0}

The diagram commutes by the equality of the two maps
$B_i\arrow A_1\otimes B_{i-2}\otimes A_1$ mentioned above.
View this diagram as having
the group $A_0\otimes B_0\otimes A_{j+l}$ in position $(0,0)$.
Multiplying the vertical maps in odd columns by $-1$ makes
this commutative diagram into a double complex $C$, meaning that
the two composite maps in each square add up to zero. Compare
the two spectral sequences converging to the cohomology of the total
complex $\Tot(C)$ \cite[section 2.4]{McCleary}. In the first one,
${}_I\! E_0^{pq}=C^{pq}$ and the differential $d_0:{}_I\! E_0^{pq}
\arrow {}_I\! E_0^{p,q+1}$ is the vertical
differential of $C$. Column $p$ of $C$
(for $0\leq p\leq l$) is $A_p$ tensored with
the Koszul complex (1) in degree $j+l-p$,
truncated at the step $\min\{ l-p,N\}$. Therefore,
$${}_I\! E_1^{pq}= \begin{cases} A_p \otimes \ker(B_{l-p}\otimes A_j
\arrow B_{l-p+1}\otimes A_{j+1}) &\text{if }q=-l\text{ and }
p+q\geq -(N-1)\\
0 &\text{if }q\neq -l\text{ and }p+q\geq -(N-1)\\
? &\text{if }p+q\leq -N.
\end{cases}$$
So $H^p(\Tot(C))$ is isomorphic to the cohomology of the complex
$$A_{l-N+1}\otimes \ker(B_{N-1}\otimes A_j\arrow B_{N-2}\otimes
A_{j+1})\arrow \cdots\arrow A_{l-1}\otimes \ker(B_1\otimes A_j
\arrow B_0\otimes A_{j+1})\arrow A_l\otimes A_j\arrow 0$$
for $-(N-2)\leq p\leq 0$, where $A_l\otimes A_j$ is placed in degree zero.
(In this range, these groups in the $E_1$ term
cannot be hit by any differential after $d_1$.)

The second spectral sequence converging to $H^*(\Tot(C))$
has ${}_{II}\! E_0^{pq}=C^{qp}$, and the differential $d_0:
{}_{II}\! E_0^{pq}\arrow {}_{II}\! E_0^{p,q+1}$ 
corresponds to the horizontal differential
in $C$. Row $-r$ in $C$ is the Koszul complex (1')
in degree $r$ (with the group
$k$ omitted in the case $r=0$),
truncated at the $N$th step, and tensored with
$A_{j+l-r}$. So 
$${}_{II}\! E_1^{pq}=\begin{cases} A_{j+l} &\text{if }p=q=0\\
0 &\text{if }p+q\geq -(N-1)\text{ and }(p,q)\neq (0,0)\\
? &\text{if }p+q\leq -N.
\end{cases}$$
Therefore $H^p(\Tot(C))$ is isomorphic to $A_{j+l}$ if $p=0$
and to zero if $-(N-1)\leq p\leq -1$ (although we only need this
for $-(N-2)\leq p\leq -1$). Combining this with the previous
description of $H^p(\Tot(C))$ gives the exact sequence (4),
as we want.
\qed

\section{Castelnuovo-Mumford regularity of a tensor product}

The properties of Castelnuovo-Mumford regularity discussed
in section \ref{notation}
follow immediately from the classical case of sheaves
on projective space. A deeper fact is that,
{\it since $O_X(1)$ is Koszul-ample, }regularity behaves well
under tensor products (Theorem \ref{product}). This will be used
along with Lemma \ref{rvan} in the proof of our main
vanishing theorem, Theorem \ref{van}. The theorem that regularity
behaves well under tensor products was proved
by Arapura \cite[section 1]{ArapuraPart}, but the statements there
have be corrected slightly (we have to define Koszul-ampleness
to be in degrees out to $2n$, not just $n$,
for the proof of \cite[Lemma 1.7]{ArapuraPart} to work).
Also, Theorem \ref{orlov}
simplifies the Koszulity assumption needed for these
results, and \cite[section 1]{ArapuraPart}
works with a smooth variety, an assumption which can be dropped.
So it seems reasonable to give the proofs here.

Throughout this section, let $X$ be a projective scheme of dimension
$n$ over a field such that the ring $O(X)$ is a field (for example,
any connected reduced projective 
scheme over a field). Write $k=O(X)$. Let $O_X(1)$
be a very ample line bundle, and
define the vector bundles $\RR_i$ as in section \ref{koszul}. (We will
only consider $\RR_i$ when $O_X(1)$ is at least $i$-Koszul.)

\begin{lemma}
\label{divide}
Let $E$ be a vector bundle and $F$ a coherent sheaf on $X$.
Let $i\geq 0$, and assume that $O_X(1)$ is $(2n-i+1)$-Koszul. Suppose
that for each pair of integers $0\leq a\leq 2n-i$ and $b\geq 0$,
either $H^b(X,E\otimes \RR_a)=0$ or $H^{i+a-b}(X,F(-a))=0$. Then
$$H^i(X,E\otimes F)=0.$$
\end{lemma}

{\bf Proof. }This is essentially \cite[Lemma 1.6]{ArapuraPart}. 
Theorem \ref{orlov} gives the first $2n-i$ steps of 
a resolution of the diagonal on $X\times_k X$. Tensoring with $E\boxtimes F$
gives a resolution of $E\otimes F$ on the diagonal in $X\times X$:
$$(E\otimes \RR_{2n-i})\boxtimes F(-2n+i)\arrow \cdots\arrow
(E\otimes \RR_0)\boxtimes F(0)\arrow E\otimes F\arrow 0.$$
To check that the latter complex really is exact, we have to show
that the sheaves $\Tor_i^{O_X\otimes_k O_X}(E\otimes_k F, O_X)$
are zero for $i>0$. Since $E$ is a vector bundle, it suffices to show
that $\Tor_i^{O_X\otimes_k O_X}(O_X\otimes_k F, O_X)=0$ for $i>0$.
But this is isomorphic to $\Tor_i^{O_X}(F, O_X)=0$, which is indeed
zero for $i>0$.

It follows that $H^i(X,E\otimes F)$ is zero if we have
$H^{i+a}(X\times X,(E\otimes \RR_a)\boxtimes F(-a))=0$
for $0\leq a\leq 2n-i$. (Because $X\times X$ has dimension $2n$,
it does not matter how the resolution continues beyond
degree $2n-i$.) This vanishing follows from our assumption
by the K\"unneth formula. \qed

For $q=0$, the following is a well-known property of Castelnuovo-Mumford
regularity \cite[v.~1, Theorem 1.8.5]{Lazarsfeld}.

\begin{lemma}
\label{shift}
Let $O_X(1)$ be a basepoint-free line bundle on a projective
scheme $X$ over a field. Let $F$ be a coherent sheaf and $q$ a natural number
such that
\begin{equation*}
0=H^{q+1}(X,F(-1))=H^{q+2}(X,F(-2))=\cdots
\tag{$C_q$}
\end{equation*}
Then $F(1)$ also satisfies $C_q$. That is,
$0=H^{q+1}(X,F)=H^{q+2}(X,F(-1))=\cdots$.
\end{lemma}

{\bf Proof. }Let $B_1=H^0(X,O(1))$, and let $\RR_1$ be the kernel:
$$0\arrow \RR_1\arrow B_1\otimes O_X\arrow O_X(1)\arrow 0.$$
Clearly $\RR_1$ is a vector bundle. We prove the lemma by
descending induction on $q$, it being trivial for $q\geq n=\dim(X)$.
So suppose the result holds with $q$ replaced by $q+1$. For each integer
$j$, we have an exact sequence $0\arrow \RR_1\otimes F(-j-1)
\arrow B_1\otimes F(-j-1)\arrow F(-j)\arrow 0$ of sheaves on $X$.
This gives a long
exact sequence of cohomology, for each $j\geq 1$:
$$H^{q+j}(X,F(-j))\arrow H^{q+j+1}(X,\RR_1\otimes F(-j-1))
\arrow B_1\otimes H^{q+j+1}(X,F(-j-1)).$$
The groups on the left and right are zero by our assumption $C_q$,
and so $H^{q+1+j}(X,\RR_1\otimes F(-j-1))=0$ for all $j\geq 1$.
That is, the sheaf $\RR_1\otimes F(-1)$ satisfies condition
$C_{q+1}$. By our inductive hypothesis, it follows
that $\RR_1\otimes F$ also satisfies $C_{q+1}$.

Now use a different
twist of the short exact sequence of sheaves above to get:
$$B_1\otimes H^{q+j}(X,F(-j))\arrow H^{q+j}(X,F(-j+1))\arrow
H^{q+j+1}(X,\RR_1\otimes F(-j)).$$
For all $j\geq 1$, the left group is zero by assumption and
the right group is zero as we have just shown. So the middle
group is zero, for all $j\geq 1$. That is, $F(1)$ satisfies
condition $C_q$. \qed

We return to our standing assumptions in this section:
$O_X(1)$ is a very ample
line bundle on a projective scheme $X$ such that the ring
$O(X)$ is a field $k$,
and $X$ has dimension $n$ over $k$.

\begin{lemma}
\label{rvan}
Let $F$ be a coherent sheaf on $X$. Let $q$ be a natural number such
that
$$H^{q+1}(X,F(-1))=H^{q+2}(X,F(-2))=\cdots.$$
Then $H^j(X,\RR_i\otimes F)=0$ for $j>q$
if $O_X(1)$ is $(n-j+1+i)$-Koszul.
\end{lemma}

{\bf Proof. }This is \cite[Corollary 1.9]{ArapuraPart}, with the
Koszulity assumption corrected. The method is to show more generally
that for any $a\geq 0$, $i\geq 0$, and $q+a<j$, $H^j(X,\RR_i\otimes
F(-a))=0$ if $O_X(1)$ is $(n-j+1+i)$-Koszul.
We prove this by descending
induction on $a$, starting with $a\geq n$ where the result is automatic
(since $j>q+a\geq n$).
We can assume that $j\leq n$, otherwise the cohomology
group we consider is automatically zero. Therefore,
$N:=n-j+1+i$ is at least $i+1$.
The sequence (2) in section \ref{koszul} gives an exact sequence
$$0\arrow \RR_{i+1}\otimes O_X(-1)\arrow B_{i+1}\otimes O_X\arrow
\RR_i\arrow 0,$$
since $i+1\leq N$. This gives an exact sequence of cohomology
$$B_{i+1}\otimes H^j(X,F(-a-1))\arrow H^j(X,\RR_i\otimes F(-a))
\arrow H^{j+1}(X,\RR_{i+1}\otimes F(-a-1)).$$
The first group is zero by our assumption on $F$ together with
Lemma \ref{shift}, and the last group is zero by our descending
induction on $a$. Thus $H^j(X,\RR_i\otimes F(-a))=0$ as we want.
\qed

We now generalize a standard property of Castelnuovo-Mumford
regularity from sheaves on projective space to sheaves on an arbitrary
reduced projective scheme. We follow
Arapura's proof \cite[Corollary 1.12]{ArapuraPart}, with the
Koszulity assumption corrected.

\begin{theorem}
\label{product}
Let $X$ be a projective scheme of dimension $n$ over a field
such that the ring of regular functions on $X$ is a field 
(example: $X$ connected and reduced).
Let $O_X(1)$ be a $2n$-Koszul line bundle on $X$. Let $E$ be a vector
bundle and $F$ a coherent sheaf on $X$. Then
$$\reg(E\otimes F)\leq \reg(E)+\reg(F).$$
\end{theorem}

{\bf Proof. } By definition, $E(\reg(E))$ and
$F(\reg(F))$ have regularity zero. So we can assume that $E$ and $F$
have regularity at most zero, and we want to show that $E\otimes F$
has regularity at most zero. That is, we have to show that
$$H^{i}(X,E\otimes F(-i))=0$$
for $i>0$. For any $0\leq a\leq n+b-1$ and $b>0$,
we have $H^b(X,E\otimes \RR_a)=0$
by Lemma \ref{rvan} (in particular, we have arranged for the Koszulity
assumption in Lemma \ref{rvan} to hold).
For any $b>0$ and any $n+b\leq a\leq 2n-i$, we have
$$H^{i+a-b}(X,F(-i-a))=0$$
for dimension reasons. Finally, for any $0\leq a\leq 2n-i$ and $b=0$,
$$H^{i+a-b}(X,F(-i-a))=0$$
since $F$ has regularity at most zero. 
Then Lemma \ref{divide} gives that
$H^i(X,E\otimes F(-i))=0$ for $i>0$. \qed

\section{Hochschild homology and the Frobenius homomorphism}

In this section we prove a flatness property of the Frobenius
homomorphism $F(x)=x^p$
for all commutative $\F_p$-algebras, to be used in the proof
of our main vanishing theorem, Theorem \ref{van}. This seems striking,
since the Frobenius homomorphism on a noetherian $\F_p$-algebra $A$
is flat only for $A$ regular \cite[Corollary 2.7]{Kunz}.
The statement is that the $\Tor$ groups of certain $A\otimes A$-modules
vanish in positive degrees,
or equivalently that certain Hochschild homology groups vanish
in positive degrees. The proof
is inspired by Pirashvili's proof of vanishing for a related set 
of $\Tor$ groups in positive characteristic \cite{Pirashvili}.

Let $k$ be a field of characteristic $p>0$, and let $N$
be a natural number. For a commutative $k$-algebra $A$, write
$\widetilde{A}=k\otimes_k A$, where $k$ maps to $k$ by the $N$th iterate of 
Frobenius, $x\mapsto x^{p^N}$. We view $\widetilde{A}$
as a $k$-algebra
using the left copy of $k$. The {\it relative
Frobenius homomorphism }$\varphi:\widetilde{A}\arrow A$
is the $k$-algebra homomorphism given by $\varphi(x\otimes y)=
x\otimes y^{p^N}$ for $x\in k$, $y\in A$.

\begin{theorem}
\label{flat}
For any commutative $k$-algebra $A$, 
view $\widetilde{A}\otimes_k A$ as a module
over $\widetilde{A}\otimes_k \widetilde{A}$ by
$x\otimes y\mapsto x\otimes \varphi(y)$,
and $\widetilde{A}$ as a module
over $\widetilde{A}\otimes_k \widetilde{A}$ by $x\otimes y\mapsto xy$. Then
$$\Tor_i^{\widetilde{A}\otimes_k \widetilde{A}}(\widetilde{A}\otimes_k A, \widetilde{A})\cong
\begin{cases}
A & \text{if }i=0\\
0 & \text{if }i>0.
\end{cases}$$
\end{theorem}

The Hochschild homology of a $k$-algebra $R$ with coefficients in an
$R$-bimodule $M$ can be defined as $H_i(R,M)=\Tor_i^{R\otimes_k R^{\op}}
(M,R)$ \cite[Proposition 1.1.13]{Loday}.
So an equivalent formulation of Theorem \ref{flat} is that
$$H_i(\widetilde{A},\widetilde{A}\otimes_k A)\cong
\begin{cases}
A & \text{if }i=0\\
0 & \text{if }i>0.
\end{cases}$$

{\bf Proof. }In this proof,
all tensor products are over $k$ unless otherwise specified.

The theorem is easy in degree zero: an isomorphism
$$(\widetilde{A}\otimes A)\otimes_{\widetilde{A}\otimes \widetilde{A}} \widetilde{A}\arrow A$$
is given by mapping $x\otimes y$ to $\varphi(x)\otimes y$. 

The vanishing of $\Tor$ in positive degrees
is clear for $R$ a free commutative
$k$-algebra (that is, a polynomial ring, possibly on infinitely many
variables). Indeed, the relative
Frobenius homomorphism $\varphi:\widetilde{A}\arrow A$ is flat in that case,
so that $\widetilde{A}\otimes A$ is a flat $\widetilde{A}\otimes \widetilde{A}$-module.

For an arbitrary commutative $k$-algebra, let $P_*\arrow A$ be a
{\it free resolution }of $A$, meaning an $A$-augmented simplicial
commutative $k$-algebra $P_*\arrow A$ which is acyclic and such
that each $P_i$ is a free commutative $k$-algebra. This exists
\cite[section 3.5.1]{Loday}.

We are trying to show that Hochschild homology $H_i(\widetilde{A},\widetilde{A}\otimes_k A)$
is zero for $i>0$. The standard complex computing Hochschild homology
$H_i(R,M)$ for a $k$-algebra $R$ and an $R$-bimodule $M$
consists of the $k$-vector spaces $C_n(R,M)=M\otimes R^{\otimes n}$.
These form a simplicial $k$-vector space, with boundary
maps $d_0:M\otimes R^{\otimes n}\arrow M\otimes R^{\otimes n-1}$
given by \cite[section 1.1.1]{Loday}:
\begin{align*}
d_0(m,a_1,\ldots,a_n)&=(ma_1,a_2,\ldots,a_n)\\
d_i(m,a_1,\ldots,a_n)&=(m,a_1,\ldots,a_ia_{i+1},\ldots,a_n)
\text{ for }1\leq i\leq n-1,\\
d_n(m,a_1,\ldots,a_n)&=(a_nm,a_1,\ldots,a_{n-1}).
\end{align*}
So we can view Hochschild homology $H_i(R,M)$ as the homology
groups of the chain complex $C_*(R,M)$ with boundary map
$$b=\sum_{i=0}^n (-1)^id_i.$$

For any homomorphism $R\arrow S$ of $k$-algebras and homomorphism
$M\arrow N|_R$ of $R$-bimodules, we have an obvious homomorphism
$C_*(R,M)\arrow C_*(S,N)$ of the Hochschild chain complexes.
Imitating the proof of \cite[Theorem 3.5.8]{Loday},
we consider the homomorphisms between
the chain complexes $C_*(P_n,\widetilde{P_n}\otimes P_n)$
given by the
$k$-algebra homomorphisms $d_i:P_n\arrow P_{n-1}$ for $0\leq i\leq n$.
We get a commutative diagram, where the horizontal arrows
are the alternating sum of the homomorphisms given by $d_0,\ldots,d_n$:
$$\begin{CD}
@VVV  @VVV  @VVV  \\
(\widetilde{P_0}\otimes P_0)\otimes \widetilde{P_0}^{\otimes 2} 
@<<< (\widetilde{P_1}\otimes P_1)\otimes \widetilde{P_1}^{\otimes 2} 
@<<< (\widetilde{P_2}\otimes P_2)\otimes \widetilde{P_2}^{\otimes 2} 
@<<< \\
@VbVV  @VbVV  @VbVV  \\
(\widetilde{P_0}\otimes P_0)\otimes \widetilde{P_0} 
@<<< (\widetilde{P_1}\otimes P_1)\otimes \widetilde{P_1} 
@<<< (\widetilde{P_2}\otimes P_2)\otimes \widetilde{P_2} 
@<<< \\
@VbVV  @VbVV  @VbVV  \\
\widetilde{P_0}\otimes P_0
@<<< \widetilde{P_1}\otimes P_1
@<<< \widetilde{P_2}\otimes P_2
@<<<
\end{CD}$$
Define a double complex $M$
by changing $b$ to $-b$ in the odd columns of this
commutative diagram.

The homology of the $i$th column in the double complex $M$ is
the Hochschild homology of the polynomial ring $P_i$ which
we have already computed:
$$H_j(\widetilde{P_i},\widetilde{P_i}\otimes P_i)=\begin{cases}
P_i &\text{if }j=0\\
0 & \text{if }j>0.
\end{cases}$$
Therefore the homology of the total complex, $H_i(\Tot(M))$,
is isomorphic to
$$H_i(P_0\leftarrow P_1 \leftarrow P_2 \leftarrow \cdots)
\cong \begin{cases} 
A & \text{if }i=0\\
0 & \text{if }i>0.
\end{cases}$$

On the other hand, the homology of the $j$th row of the double complex
$M$ is
$$H_i(\widetilde{P_*}\otimes P_*\otimes \widetilde{P_*}^{\otimes j})
\cong \begin{cases}
\widetilde{A}\otimes A\otimes \widetilde{A}^{\otimes j} &\text{if }i=0\\
0 &\text{if }i>0.
\end{cases}$$
This calculation uses the Eilenberg-Zilber theorem: when we define the tensor
product of two simplicial $k$-vector spaces by $(A\otimes B)_n=A_n
\otimes_k B_n$, with boundary maps $d_i(x\otimes y)=d_i(x)
\otimes d_i(y)$, the associated chain complex has homology groups
$H_i(A\otimes B)$ isomorphic to
$H_i(A)\otimes H_i(B)$ \cite[Theorem 1.6.12]{Loday}.
By the second spectral sequence converging to the homology
of $\Tot(M)$, we have
an isomorphism
$$H_i(\widetilde{A}\otimes A\leftarrow (\widetilde{A}\otimes A)
\otimes \widetilde{A}\leftarrow \cdots )\cong \begin{cases}
A & \text{if }i=0\\
0 & \text{if }i>0.
\end{cases}$$
The complex here (a quotient of the 0th column of $M$)
is the one that computes Hochschild homology $H_i(\widetilde{A},
\widetilde{A}\otimes A)$. Thus we have shown that
$$H_i(\widetilde{A},\widetilde{A}\otimes A)\cong \begin{cases}
A & \text{if }i=0\\
0 & \text{if }i>0.
\end{cases}$$
\qed

\section{Vanishing in positive characteristic}

In this section, we prove a vanishing theorem for
reduced projective schemes over a field of positive characteristic. For smooth
projective varieties, the theorem is due to Arapura
\cite[Theorem 5.4]{ArapuraPart}.
The generalization to singular schemes follows the original
proof, but with a new ingredient, a flatness property
of the Frobenius morphism for arbitrary schemes over $\F_p$
(Theorem \ref{flat}).

Theorem \ref{van} is very different from the best-known
vanishing theorem in positive characteristic, Deligne-Illusie-Raynaud's
version of the Kodaira vanishing theorem.
Their proof shows that
(when a smooth projective variety $X$ lifts from $\Z/p$ to $\Z/p^2$)
vanishing of cohomology
groups for the line bundle $K_X\otimes L^{p^b}$ with $b$ large
can imply vanishing for $K_X\otimes L$ \cite[Lemme 2.9]{DI}.
Theorem \ref{van} goes the opposite way.

One result related to Theorem \ref{van}
is Siu's nonvanishing theorem.
Siu's theorem says that if $E$ is a pseudoeffective line bundle on a smooth
projective variety $X$ of dimension $n$ over a field of characteristic
zero, and $O_X(1)$ is an ample line bundle on $X$, then
$H^0(X,K_X\otimes E(j))\neq 0$ for some $1\leq j\leq n+1$. 
(Ein pointed out this slight extension of Siu's theorem as stated
in \cite[Corollary 9.4.24]{Lazarsfeld}.)
By Serre duality, Siu's theorem gives that
if $L$ is a line bundle with $0=H^n(X,L(-1))
=\cdots =H^n(X,L(-n-1))$, then $L$ is naively $(n-1)$-ample (that is,
$L^*$ is not pseudoeffective, by Theorem \ref{n-1}).

We use the notion
of regularity $\reg(M)$ from section \ref{notation}.

\begin{theorem}
\label{van}
Let $X$ be a projective
scheme of dimension $n$
over a field of characteristic $p>0$ such that $O(X)$ is a field
(example: $X$ connected and reduced). Let $O_X(1)$
be a Koszul-ample
line bundle on $X$. Let $q$ be a natural number.
Let $L$ be a line bundle on $X$ with
$$0=H^{q+1}(X,L(-n-1))=H^{q+2}(X,L(-n-2))=\cdots .$$
Then for any coherent sheaf $M$ on $X$, we have
$$H^i(X,L^{\otimes p^b}\otimes M)=0$$
whenever $i>q$ and $p^b\geq \reg(M)$.
\end{theorem}

{\bf Proof. }We follow Arapura's proof in the smooth case as far as possible.
Let $k$ be the field $H^0(X,O_X)$.
Let $i$ be an integer greater than $q$ and $b$ a natural number
such that $p^b\geq \reg(M)$.
Let $f'=F^b$, where $F:X\arrow X$ is the absolute Frobenius
morphism (which acts as the identity on $X$ as a set, and acts by
$p$th powers on $O_X$). Let $\widetilde{X}$ be the base extension $X\times_k k$
where $k$ maps to $k$ by $x\mapsto x^{p^b}$. Then $f'$ factors as
$$X\stackrel{f}{\longrightarrow} \widetilde{X}
  \stackrel{g}{\longrightarrow} X,$$
where $f$ is a morphism of $k$-schemes and $g$ is the natural morphism.
The morphism $f$ is called the relative $b$th Frobenius morphism. Let
$\widetilde{L}=g^*L$ and $O_{\widetilde{X}}(1)=g^*O_X(1)$. Since $g$ is
given by a field extension, $O_{\widetilde{X}}(1)$ is Koszul-ample.

Let
$$C^i=\begin{cases} 
\widetilde{\RR}_{-i}
\boxtimes O_{\widetilde{X}}(i) & \text{if }-2n+1\leq i\leq 0 \\
\ker(\widetilde{\RR}_{2n-1}\boxtimes O_{\widetilde{X}}(-2n+1)\arrow 
\widetilde{\RR}_{2n-2}
\boxtimes O_{\widetilde{X}}(-2n+2) & \text{if }i=-2n,
\end{cases}$$
where $\RR_i$ is defined as in section \ref{koszul}
and $\widetilde{\RR}_i=g^*\RR_i$. The sheaves $C^i$
form a resolution $C^*$ of the diagonal $\Delta$ on $\widetilde{X}
\times \widetilde{X}$, by Theorem \ref{orlov}.
More generally, for vector bundles $E_1$
and $E_2$ on $\widetilde{X}$, $(E_1\boxtimes E_2)\otimes C^*$
is quasi-isomorphic to $\delta_*(E_1\otimes E_2)$, where
$\delta:\X\arrow \X\times \X $ denotes the diagonal embedding.
Therefore, $D^*=(O_{\X}(-n)\boxtimes O_{\X}(n))\otimes C^*$
is another resolution of the diagonal.

Let $\Gamma\subset \X\times X$ be the graph of $f:X\arrow \X$,
and write $\gamma$ for the inclusion $X\cong \Gamma\subset \X\times X$.
Then $\Gamma=(1\times f)^{-1}(\Delta)$. Since $D^*$ is a resolution
of $O_{\Delta}\cong O_{\X}$, the cohomology sheaves of $(1\times f)^*(D^*)$
are the groups $\Tor_*^{O_{\X\times \X}}(O_{\X\times X},O_{\X})$.
Since $X$ may be singular, the relative Frobenius morphism
need not be flat, and so $O_{\X\times \X}$ need not be flat over
$O_{\X\times X}$. Nonetheless, Theorem \ref{flat} shows
that these Tor groups are zero in positive degrees. Therefore,
$(1\times f)^*(D^*)$ is a resolution of $O_{\Gamma}$ on $\X\times X$.

I claim that the complex of sheaves
$$G^*=(\widetilde{L}\boxtimes M)\otimes (1\times f)^*(D^*)$$
is a resolution of the sheaf
$\gamma_*(f^*\widetilde{L}\otimes M) \cong \gamma_*(L^{\otimes p^b}\otimes M)$.
This follows if we can show that the sheaf
$\Tor_i^{O_{\X}\otimes_k O_X}(\widetilde{L}\otimes_k M, O_{\Gamma})$ is zero
for $i>0$, where $O_{\Gamma}\cong O_X$. Since $\widetilde{L}$
is a line bundle, it suffices to show that 
$\Tor_i^{O_{\X}\otimes_k O_X}(O_{\X}\otimes_k M, O_{X})$ is zero
for $i>0$. But this is isomorphic to 
$\Tor_i^{O_{X}}(M, O_{X})$, which is indeed zero for $i>0$.
Therefore, we can compute $H^i(X,L^{\otimes p^b}\otimes M)=H^i(\X\times X,
\gamma_*(L^{\otimes p^b}\otimes M))$ using the resolution $G^*$.

By the spectral sequence $E_1^{i+c,-c}=H^{i+c}(\X\times X,G^{-c})\imp H^{i}
(\X\times X,G^*)$, the theorem holds if we can show
that $H^{i+c}(\X\times X,G^{-c})=0$ for $i>q$ and $c\geq 0$. This is clear
for the leftmost sheaf in the resolution, corresponding to $c=2n$,
because $\X\times X$
has dimension only $2n$ and $q$ is nonnegative. For $0\leq c\leq 2n-1$,
the K\"unneth formula gives that
\begin{align*}
H^{i+c}(\X\times X,G^{-c})&= H^{i+c}(\X\times X,
(\widetilde{\RR}_c\otimes\widetilde{L}(-n))\boxtimes
M(p^b(n-c)))\\
&\cong \oplus_{r+s=i+c} H^r(\X,\widetilde{\RR}_c\otimes\widetilde{L}(-n))
\otimes H^s(X,M(p^b(n-c))).
\end{align*}
It remains to show that for all $i>q$, all $0\leq c\leq 2n-1$,
and all $r+s=i+c$, either the $H^r$ group or the $H^s$ group here is zero.

First suppose that $r>q$. We can assume that $s\leq n$; otherwise
the $H^s$ group is zero. So $r\geq i+c-n$, or equivalently
$c\leq n+r-i$, and so $c\leq n+r-1$. By that fact together with
$r>q$, Koszul-ampleness of $O_X(1)$, and our assumption on $L$,
Lemma \ref{rvan} gives that $H^r(X,\RR_c\otimes L(-n))=0$.
The implies the same vanishing on $\X$ for the same sheaf pulled back
by the base extension $\X\arrow X$.

Otherwise, $r\leq q$. Since $r+s=i+c$ and $i>q$, we have $s>c$.
In particular, $s$ is greater than zero.
We can assume that $s\leq n$; otherwise
the $H^s$ group is zero. So $c<s\leq n$. Then
$p^b(n-c)\geq p^b\geq \reg(M)$, and so (since $s$ is greater than zero)
$H^s(X,M(p^b(n-c)))=0$. \qed

\section{$q$-T-ampleness}

{\bf Definition. }
Let $O_X(1)$ be a Koszul-ample line bundle (defined in section
\ref{notation}) on a projective scheme
$X$ of dimension $n$ over a field such that $O(X)$ is a field
(example: $X$ connected and reduced).
Let $q$ be a natural number.
A line bundle $L$ on $X$ is called
{\it $q$-T-ample }if, for some positive integer $N$, we have
$$0=H^{q+1}(X,L^{\otimes N}(-n-1))=H^{q+2}(X,L^{\otimes N}(-n-2))=\cdots
= H^n(X,L^{\otimes N}(-2n+q)).$$

The details are not too important.
Given that
some multiple of $L$ kills cohomology above dimension $q$ for
finitely many explicit line bundles on $X$, as in this definition,
we will deduce that infinitely many
multiples of $L$ kill cohomology above dimension $q$ for any given coherent
sheaf. In particular, that will show that
the definition of $q$-T-ampleness of a line bundle
is independent of the choice of the Koszul-ample line bundle $O_X(1)$.

Theorem \ref{van}
makes it easy to characterize $q$-T-ampleness in characteristic
$p>0$ by a certain asymptotic vanishing of cohomology, as follows.
In particular, property (2) shows that $q$-T-ampleness in characteristic $p>0$
does not depend on the choice of Koszul-ample line bundle $O_X(1)$.

\begin{corollary}
\label{charp}
Let $X$ be a projective scheme over a field of characteristic $p>0$
such that $O(X)$ is a field (example: $X$ connected and reduced).
Let $q$ be a natural number. The following properties are equivalent,
for a line bundle $L$ on $X$.

(1) $L$ is $q$-T-ample.

(2) Some positive multiple of $L$ is $q$-F-ample. That is, there is
a positive integer $N$ such that for all coherent sheaves $M$ on $X$,
we have $H^i(X,M\otimes L^{\otimes Np^j})=0$ for all $i>q$ and 
all $j$ sufficiently large depending on $M$.
\end{corollary}

{\bf Proof. }It is immediate that (2) implies (1). Here the number $N$
in the definition of $q$-T-ampleness
will be of the form $Np^j$ for some $j$ in the notation of (2). 
Theorem \ref{van} shows that (1) implies (2), with the same value
of $N$. \qed

In characteristic zero, $q$-T-ampleness has even stronger consequences,
as we now show. Curiously, the proof involves reducing modulo many different
prime numbers and using the prime number theorem.
In particular, Theorem \ref{equivred} shows that 
$q$-T-ampleness in any characteristic
is independent of the choice of the Koszul-ample line bundle $O_X(1)$.

\begin{theorem}
\label{equivred}
Let $X$ be a projective scheme over a field of characteristic zero
such that $O(X)$ is a field (example: $X$ connected and reduced).
Let $q$ be a natural number. The following properties are equivalent,
for a line bundle $L$ on $X$.

(1) $L$ is $q$-T-ample.

(2) $L$ is naively $q$-ample. That is, for every coherent sheaf
$M$ on $X$, we have $H^i(X,M\otimes L^{\otimes m})=0$ for
all $i>q$ and all $m$ sufficiently large depending on $M$.

(3) $L$ is uniformly $q$-ample.
That is, there is a constant $\lambda>0$ such that
$H^i(X,L^{\otimes m}(-j))=0$ for all $i>q$, $j>0$, and $m\geq \lambda
j$.
\end{theorem}

{\bf Proof. }Demailly-Peternell-Schneider showed that (3) implies (2),
by resolving any coherent sheaf by direct sums of the line bundles
$O_X(-j)$ \cite[Proposition 1.2]{DPS}. Clearly (2) implies (1).
We will show that (1) implies (3).
Let $L$ be a $q$-T-ample line bundle, and let $N$
be the positive integer given in the definition.

To prove (3), we can work on a model of $X$ over some finitely
generated field extension of $\Q$. We can extend this to a projective model
$X_R$
of $X$ over some domain $R$ which is a
finitely generated $\Z$-algebra. We can assume that $R=O(X_R)$, after
replacing $R$ by a finite extension if necessary;
then all fibers $X_t$ over closed points $t$ of $\Spec(R)$
have $O(X_t)$ equal to 
a field (the residue field at $t$). After inverting a nonzero element of $R$,
we can assume that $X_R$ is also flat
over $R$
and $O_X(1)$ is Koszul-ample over $R$. (Recall from section
\ref{notation} that Koszul-ampleness is a Zariski-open
condition on a line bundle.)
Choose an extension
of the line bundle $L$ to $X$ over $R$.

By semicontinuity of cohomology \cite[Theorem III.12.8]{Hartshorne},
after inverting a nonzero element of $R$, we have
$$0=H^{q+1}(X_t,L^{\otimes N}(-n-1))=\cdots = H^n(X_t,L^{\otimes N}(-2n+q))$$
for all closed points $t\in\Spec(R)$, since that is true in characteristic
zero. Let $r=\max\{ 1,\reg(L)\} $ (computed in characteristic zero).
After inverting a nonzero element of $R$ again, we can assume that
$\reg(L|_{X_t})\leq r$ for all closed points $t\in \Spec(R)$, again
by semicontinuity.
By Theorem \ref{van}, for
each line bundle $M$ on $X_R$ and each closed point
$t\in\Spec(R)$, we have
$$H^i(X_t,L^{\otimes Np^b}\otimes M|_{X_t})=0$$
whenever $i>q$ and $p^b\geq \reg(M|_{X_t})$,
where $p$ denotes the characteristic of the (finite) residue field
of $R$ at $t$.
Again by semicontinuity of cohomology,
it follows that for each closed point $t\in\Spec(R)$ and each
line bundle $M$ on $X_R$,
we have
\begin{equation*}
H^i(X,L^{\otimes Np^b}\otimes M)=0
\tag{1}
\end{equation*}
in characteristic zero whenever $i>q$ and $p^b\geq \reg(M|_{X_t})$.

This already leads to an equivalent characterization
of $q$-T-ampleness in characteristic zero; for example, it implies
that $q$-T-ampleness is independent of the choice of Koszul-ample
line bundle $O_X(1)$. But we want to prove the stronger statements
that $L$ is naively $q$-ample and in fact uniformly $q$-ample.

The scheme $\Spec(R)$ has closed points of every characteristic $p$ at least
equal to some positive integer $p_0$.
Therefore,
for any positive integer $a$, we have
\begin{equation*}
H^i(X,L^{\otimes Np+a}(-j))=0
\tag{2}
\end{equation*}
in characteristic zero
if $i>q$, $p$ is a prime $\geq p_0$, and $j+ar\leq p$.
(Indeed, let $t$ be a closed point of characteristic $p$ in $\Spec(R)$.
We have  $j+a\,\reg(L|_{X_t})\leq p$ since we arranged that
$\reg(L|_{X_t})\leq r$. 
By Theorem \ref{product}, we have $\reg(A\otimes B)\leq \reg(A)+\reg(B)$
for line bundles in any characteristic, and so
$j+\reg(L^{\otimes a}|_{X_t})\leq p$. Equivalently,
$\reg(L^{\otimes a}(-j)|_{X_t})\leq p$, and then equation
(1) gives what we want.)

Let $c$ be a real number in $(0,1)$; we could take $c=1/2$
for the current proof, but Theorem \ref{lambda} will prove a stronger
estimate by taking $c$ close to zero.
By the prime number theorem \cite[Theorem 3.3.2]{Murty},
there is a positive integer $m_1$
(depending on $c$) 
such that for every integer $m\geq m_1$, there is a prime number
$p$ with 
$$\frac{m}{N+(c/r)}\leq p < \frac{m}{N}.$$
Taking $m_1$ big enough, we can also assume that
$p_0\leq m_1/(N+(c/r))$, so that the primes $p$ produced above
always have $p\geq p_0$.

Now, for each positive integer $j$, let $m_0=m_0(j)$ be the maximum
of $m_1$ and $\lceil (j/(1-c))(N+(c/r))\rceil $.
This will imply an inequality we want. First note that
$$\frac{j}{1-c}\bigg( N+(c/r)\bigg)\leq m_0,$$
and hence
$$\frac{j}{r}\bigg( N+(c/r)\bigg)
\leq m_0\frac{1-c}{r}.$$
Therefore
$$\bigg( m_0+(j/r)\bigg)
 \bigg( N+(c/r)\bigg)\leq m_0\bigg( N+(1/r)\bigg) ,$$
and hence
$$\frac{m_0+(j/r)}{N+(1/r)}\leq
 \frac{m_0}{N+(c/r)}.$$
It follows that the same inequality holds for all $m\geq m_0$ in place
of $m_0$. Combining this with the previous paragraph's result,
we find that for every $m\geq m_0$, there is a prime number $p\geq p_0$
with
$$\frac{m+(j/r)}{N+(1/r)}\leq p < \frac{m}{N}.$$
Equivalently, if we define an integer $a$ by $m=pN+a$, then
$a>0$ and $j+a\,r\leq p$. By equation (2),
we have shown that for every $m\geq m_0(j)$, we have
$$H^i(X,L^{\otimes m}(-j))=0$$
for $i>q$.

By our construction of the number $m_0(j)$, there is a positive
constant $\lambda>0$ such that $m_0(j)\leq \lambda j$ for all
$j> 0$. (Here $\lambda$ depends on the chosen constant $c\in (0,1)$;
we can take $c=1/2$, for example.) Thus we have shown that
for all $i>q$, $j > 0$,
and all $m\geq \lambda j$ we have
$$H^i(X,L^{\otimes m}(-j))=0.$$
That is, $L$ is uniformly $q$-ample. \qed

The proof shows something more precise.

\begin{theorem}
\label{lambda}
Let $X$ be a projective scheme over a field of characteristic zero
such that $O(X)$ is a field (example: $X$ connected and reduced).
Let $L$ be a $q$-T-ample (or naively $q$-ample, or uniformly $q$-ample)
line bundle on $X$. Then there
are positive numbers $m_1$ and $\lambda$ such that
$$H^i(X,E\otimes L^{\otimes m})=0$$
for all $i>q$,
all coherent sheaves $E$ on $X$, and all $m\geq \max\{ m_1,\lambda\,\reg(E)
\} $.
Moreover, for $N$ the number in the definition of $q$-T-ampleness,
we can take $\lambda$ to be any real number greater than $N$
(and some $m_1$ depending on $\lambda$).
\end{theorem}

Demailly-Peternell-Schneider proved this for some constant $\lambda$
when $L$ is uniformly $q$-ample \cite[Proposition 1.2]{DPS}.
The point here is that the same holds
for the a priori weaker notions of $q$-T-ampleness and naive $q$-ampleness
in characteristic zero, and that we have an explicit estimate for
the constant $\lambda$.

{\bf Proof. }The proof of Theorem \ref{equivred} shows that there
are positive constants $\lambda$ and $m_1$ such that for
all $i>q$, $j> 0$, and $m\geq \max\{ m_1,\lambda j\} $,
we have
$$H^i(X,L^{\otimes m}(-j))=0.$$
In terms of a constant $c\in (0,1)$
which we were free to choose, the proof shows that we can take $\lambda=
(1/(1-c))(N+(c/r))$, where $r=\max\{ 1,\reg(L)\} $. Thus,
by taking $c$ close to zero, we can make
$\lambda$ arbitrarily close to $N$.

Let $E$ be any coherent sheaf on $X$, and let $s=\reg(E)$,
$R=\max\{ 1,\reg(O_X)\} $. Then $E$ has a resolution by vector bundles
of the form
$$\cdots\arrow O(-s-2R)^{\oplus a_2}\arrow O(-s-R)^{\oplus a_1}
\arrow O(-s)^{\oplus a_0}\arrow E\arrow 0$$
\cite[Corollary 3.2]{ArapuraFrob}. For integers $i>q$ and
$m\geq \max\{ m_1,\lambda s\} $, we want to show that $H^i(X,E
\otimes L^{\otimes m})=0$. By the given resolution, this holds
if $H^{i+j}(X,L^{\otimes m}(-s-jR))=0$ for all $j\geq 0$.
By the previous paragraph, this holds if $m\geq \max\{ m_1,
\lambda(s+(n-q-1)R\} $, using that cohomology vanishes in dimensions
greater than $n$. This is enough to deduce the statement of the theorem,
after slightly increasing $\lambda$ and increasing $m_1$ as needed.
 \qed

\section{Projective schemes}

In this section we generalize the equivalence between the three notions
of $q$-ampleness to arbitrary projective schemes over a field
of characteristic zero.

\begin{theorem}
\label{equiv}
Let $X$ be a projective scheme over a field of characteristic
zero, with an ample line bundle $O_X(1)$.
Let $q$ be a natural number. Then there is a positive integer
$C$ such that: for all line bundles $L$ on $X$,
the following properties are equivalent.

(1) There is a positive integer $N$ such that $H^i(X,
L^{\otimes N}(-j))=0$ for all $i>q$ and all $1\leq j\leq C$.

(2) $L$ is naively $q$-ample. That is, for every coherent sheaf
$M$ on $X$, we have $H^i(X,M\otimes L^{\otimes m})=0$ for
all $i>q$ and all $m$ sufficiently large depending on $M$.

(3) $L$ is uniformly $q$-ample.
That is, there is a constant $\lambda>0$ such that
$H^i(X,L^{\otimes m}(-j))=0$ for all $i>q$, $j> 0$, and $m\geq \lambda
j$.
\end{theorem}

In contrast to the case of reduced schemes (Theorem \ref{equivred}),
we choose not to specify
the value of $C$ in condition (1). The proof gives
an explicit value for $C$, probably far from optimal.

{\bf Proof. }By the same arguments as in Theorem \ref{equivred},
(3) implies (2) and (2) implies (1), for any fixed choice of the
positive integer $C$. It remains to show that there is a positive
integer $C$ such that (1) implies (3), for all line bundles $L$
on $X$. The idea is to reduce to the case of a reduced scheme.

We can assume that $X$ is connected. Let $X_0$ be the underlying
reduced scheme of $X$. Then $k:=O(X_0)$ is a field.
After replacing the given ample line bundle
$O_X(1)$ by a positive multiple, we can assume that $O_X(1)$
restricts to a Koszul-ample line bundle on $X_0$ (by Backelin's theorem,
as in section \ref{notation}).

I claim that there is a positive integer $C$ (depending on $X$ and $O_X(1)$)
such that any line bundle on $X$ satisfying (1) restricts to a $q$-T-ample
line bundle on $X_0$. To see this, look at a resolution of $O_{X_0}$
as a sheaf of $O_X$-modules. Explicitly, if we define $s=\reg_X(O_{X_0})$
and $R=\max\{1,\reg_X(O_X)\}$, then $O_{X_0}$ has a resolution on $X$
of the form:
$$\cdots \arrow O_X(-s-2R)^{\oplus a_2}\arrow O_X(-s-R)^{\oplus a_1}
\arrow O_X(-s)^{\oplus a_0}\arrow O_{X_0}\arrow 0$$
\cite[Corollary 3.2]{ArapuraFrob}. Then, for any line bundle
$L$ on $X$, and any positive integer $c$, tensoring
this resolution with $L(-c)$ gives a resolution 
of $L(-c)|_{X_0}$ as a sheaf of $O_X$-modules:
$$\cdots \arrow L(-c-s-2R)^{\oplus a_2}\arrow
L(-c-s-R)^{\oplus a_1}
\arrow L(-c-s)^{\oplus a_0}\arrow L(-c)|_{X_0}
\arrow 0.$$
This gives a spectral sequence
$$E_1^{-u,v}=H^v(X,L(-c-s-uR)^{\oplus a_u})\imp
H^{v-u}(X_0,L(-c)).$$

Let $n$ be the dimension of $X$, and 
let $C=n+q+1+s+(n-q-1)R$ (which does not depend on the line bundle $L$).
Then the spectral sequence shows that for any line bundle $L$
on $X$ such that
$H^i(X,L(-j))=0$ for all $i>q$ and all $1\leq j\leq C$,
we have
$$0=H^{q+1}(X_0, L(-n-1))=H^{q+2}(X_0,L
(-n-2))=\cdots.$$

Thus, for the given value of $C$, any line bundle $L$ on $X$
which satisfies condition (1) restricts to a $q$-T-ample
line bundle on $X_0$. By Theorem \ref{equivred}, $L$
is uniformly $q$-ample on $X_0$. To complete the proof,
we have to show that $L$ is uniformly $q$-ample on $X$.

Let $I$ be the ideal sheaf $\ker(O_X\arrow O_{X_0})$. There is a natural
number $r$ such that $I^{r+1}=0$. Thus the sheaf $O_X$ has a filtration
with quotients $O_X/I,I/I^2,\ldots,I^{r}/I^{r+1}$, and these quotients
are all modules over $O_X/I=O_{X_0}$. 
By considering a resolution of the sheaves
$I^l/I^{l+1}$ (for $l=0,\ldots,r$)
by line bundles $O(-j)$ on $X_0$, the uniform $q$-ampleness of $L|_{X_0}$
implies that there
is a constant $\lambda>0$ such that
$H^i(X_0,L^{\otimes m}(-j)\otimes I^l/I^{l+1})=0$
for all $i>q$, all $0\leq l\leq r$, all $j>0$, and all $m>\lambda j$.
By our filtration of $O_X$, it follows that
$H^i(X, L^{\otimes m}(-j))=0$, for all  $i>q$, all $j>0$, and
all $m>\lambda j$.
That is, $L$ is uniformly $q$-ample on $X$. \qed

Therefore, on any projective scheme over a field
of characteristic zero, we can say ``$q$-ample'' 
to mean any of the equivalent conditions on a line bundle
in Theorem \ref{equiv}.
We mention a consequence of the proof:

\begin{corollary}
Let $X$ be a projective scheme over a field
of characteristic zero. Then a line bundle $L$ is $q$-ample on $X$
if and only if the restriction of $L$ to the underlying reduced
scheme $X_0$ is $q$-ample.
\end{corollary}

\section{Openness properties of $q$-ampleness}
\label{zariski}

In this section we check that $q$-T-ampleness is Zariski
open on families of varieties and line bundles over $\Z$.
It follows that $q$-ampleness is Zariski open in characteristic zero,
where we can use any of the three equivalent definitions
in Theorem \ref{equivred}.
Using Demailly-Peternell-Schneider's results, we also
find that $q$-ampleness defines an open cone (typically not convex)
in the N\'eron-Severi vector space $N^1(X)$ for $X$
of characteristic zero. Neither property was known
for naive $q$-ampleness. We discuss counterexamples
to these good properties in positive characteristic.

\begin{theorem}
\label{zariskiopen}
Let $\pi:X\arrow B$ be a flat projective morphism of schemes
over $\Z$. Suppose that $\pi$ has connected fibers in the sense
that $\pi_*(O_X)=O_B$. 
Let $L$ be a line bundle on $X$, and let $q$ be a natural
number.
Then the set of points $b$ of the scheme $B$ such that
$L$ is $q$-T-ample on the fiber over $b$ is Zariski open.
\end{theorem}

{\bf Proof. }This is straightforward. Suppose that $L$ is $q$-T-ample
on the fiber $X_b$
over a point $b\in B$. There is an affine open neighborhood $U$
of $b$ and a Koszul-ample line bundle $O_X(1)$ on the inverse image of
$U$. (Koszul-ampleness is a Zariski-open condition on a line bundle,
as discussed in section \ref{notation}). We can use this line bundle
$O_X(1)$ in the definition of $q$-T-ampleness.
(We have shown that the definition is independent of this choice.)
By definition of $q$-T-ampleness, there is a positive integer
$N$ with
$$0=H^{q+1}(X_b,L^{\otimes N}(-n-1))=\cdots = H^n(X_b,L^{\otimes N}(-2n+q)).$$
By semicontinuity of cohomology, the same is true over some
neighborhood of $b$. \qed

By Theorem \ref{equivred}, it follows that naive $q$-ampleness
and uniform $q$-ampleness are Zariski open conditions in characteristic zero.

Demailly-Peternell-Schneider showed that uniform $q$-ampleness
defines an open cone in the vector space $N^1(X)$ of $\R$-divisors
modulo numerical equivalence. For completeness, we give
the relevant definitions. Let $X$ be a projective
variety over a field of characteristic zero.
Define an $\R$-divisor $D$ on
$X$ to be
$q$-ample if $D$ is numerically
equivalent to a sum $cL+A$
with $L$ a $q$-ample divisor, $c$ a positive rational number,
and $A$ an ample $\R$-divisor.
(By definition, an $\R$-divisor is ample if it is a positive
linear combination of ample Cartier divisors
\cite[Definition 3.11]{Lazarsfeld}.)

\begin{theorem}
\label{opencone}
For any projective variety $X$ over a field
of characteristic zero,
$q$-ampleness for $\R$-divisors agrees with the earlier definitions
in the case of line bundles. Also, $q$-ampleness
defines an open cone in the real vector space $N^1(X)$.

Also, the sum of a $q$-ample $\R$-divisor and an $r$-ample
$\R$-divisor is $(q+r)$-ample.
\end{theorem}

Demailly-Peternell-Schneider proved Theorem \ref{opencone}
for uniform $q$-ampleness \cite{DPS}. This is equivalent to the
other notions of $q$-ampleness for $X$ projective
over a field of characteristic zero, by Theorem \ref{equiv}.

Theorem \ref{opencone} gives a simple insight into why the $q$-ample
cone need not be convex for $q>0$: the sum of two $q$-ample
divisors is typically $2q$-ample, not $q$-ample. An example
is the $(n-1)$-ample cone of a projective variety
of dimension $n$, which is equal to the
complement of the negative of the closed effective cone
by Theorem \ref{n-1}.
(Thus the $(n-1)$-ample cone
is the {\it complement }of a closed convex cone.)

In positive characteristic, the right notion of $q$-ampleness
remains to be found. In particular, naive $q$-ampleness
and uniform $q$-ampleness are not Zariski open conditions
in mixed characteristic, as one can check in the example
of the three-dimensional flag manifold $SL(3)/B$ over $\Z$ with $q=1$.
The figures show the (naive or uniform) $q$-ample cone in characteristic
zero and in any characteristic $p>0$, where the $1$-ample cone
is different from characteristic zero (and independent of $p$).
\begin{figure*}
\centering
\includegraphics[width=0.7\textwidth]{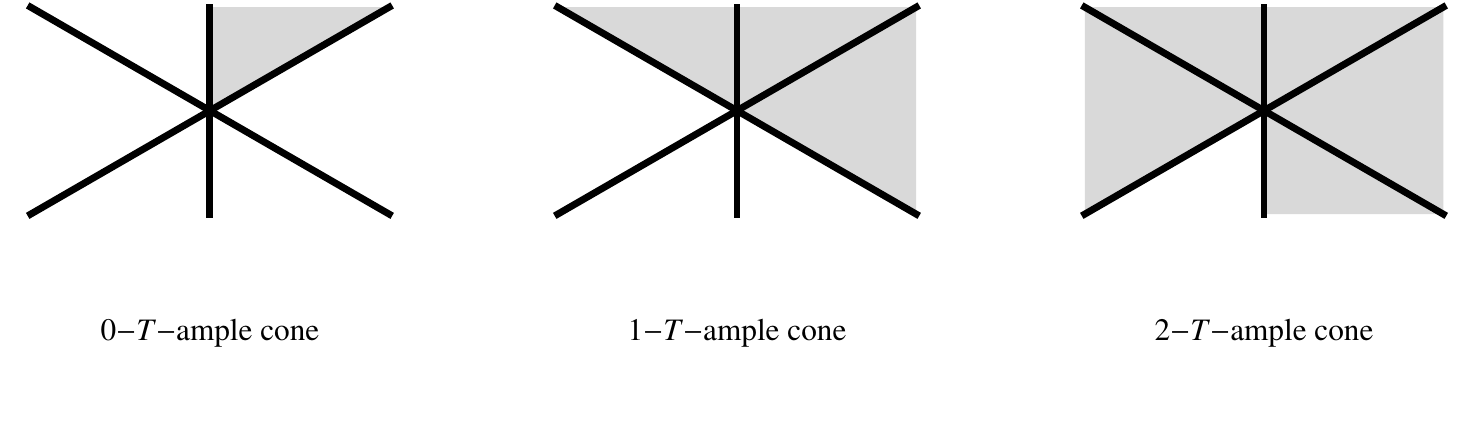}
\caption{The $q$-ample cones in $N^1(SL(3)/B)\cong \R^2$ in char.\ 0}
\end{figure*}
\begin{figure*}
\centering
\includegraphics[width=0.7\textwidth]{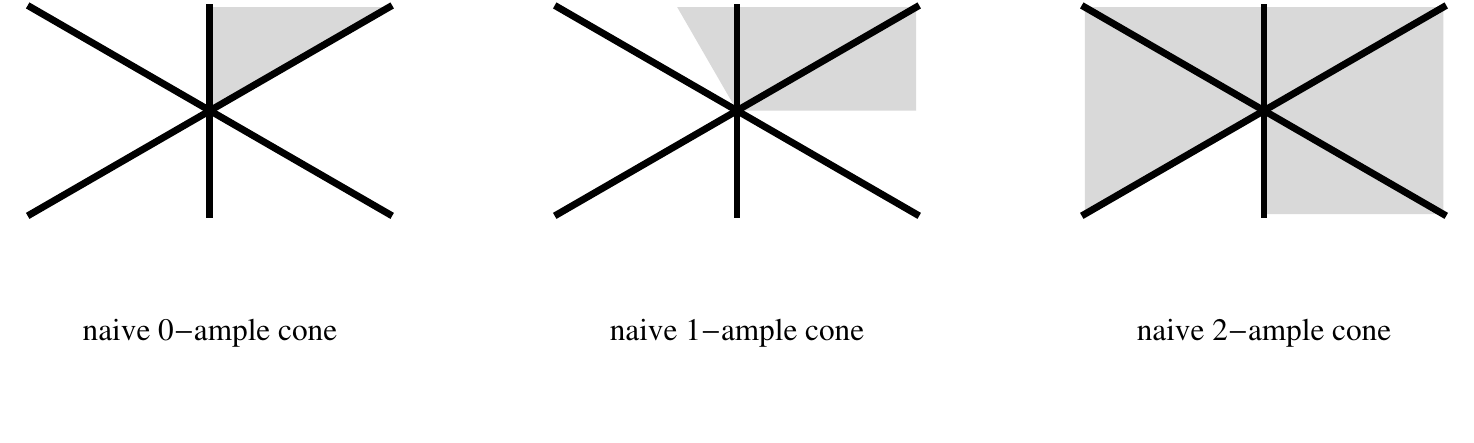}
\caption{The naive $q$-ample cones in $N^1(SL(3)/B)\cong \R^2$ in
char.\ $p>0$}
\end{figure*}

On the other hand, $q$-T-ampleness is also not a good notion in
positive characteristic, in the sense that the sum of a $q$-T-ample
divisor and an $r$-T-ample divisor need not be $(q+r)$-T-ample.
This happens with $q=r=1$ on $(SL(3)/B)^2$ in any
characteristic $p>0$.

\section{The $(n-1)$-ample cone}

The 0-ample and $(n-1)$-ample cones
of an $n$-dimensional variety are better understood than
the intermediate cases.
In this section we show that the $(n-1)$-ample cone of a projective
variety of dimension $n$ is the complement of the negative
of the closed effective cone. This is well known for smooth varieties,
but with care the proof works for arbitrary varieties.

\begin{theorem}
\label{n-1}
Let $X$ be a projective variety of dimension $n$
over a field $k$ of characteristic zero.
Let $L$ be a line bundle on $X$. Then $L$ is $(n-1)$-ample
if and only if $[L]$ in $N^1(X)$ is not in the negative of the closed
effective cone.
\end{theorem}

{\bf Proof. }Let $\omega_X$ be what Hartshorne calls the dualizing
sheaf of $X$, that is, the cohomology sheaf in dimension $n$ of the
dualizing complex of $X$. Then for every coherent sheaf $F$ on $X$,
there is a canonical isomorphism
$$\Hom_X(F,\omega_X)\cong H^n(X,F)^* $$
\cite[Proposition III.7.5]{Hartshorne}. It follows that the sheaf $\omega_X$
is torsion-free. Indeed, if $F$ is any coherent sheaf on $X$ whose
support has dimension less than $n$, then $H^n(X,F)=0$ and so
$\Hom_X(F,\omega_X)=0$.

Let $O_X(1)$ be an ample line bundle on $X$. For a coherent sheaf $F$,
write $F^*$ for the sheaf $\HHom_{O_X}(F,O_X)$. Then the sheaf $\omega_X^*$
is generically a line bundle on $X$ (in particular, it is not zero),
and so $H^0(X,\omega_X^*(j))\neq 0$ for some $j>0$. Equivalently, there
is a nonzero map $f:\omega_X\arrow O(j)$ for some $j>0$, which we fix.
Since $\omega_X$ is torsion-free, $f$ must be an injection of sheaves.

Suppose that $[L]$ in $N^1(X)$ is not in the negative of the closed
effective cone. That is, $L^*$ is not pseudoeffective. Then for any line
bundle $F$ on $X$, we have $H^0(X,F^*\otimes (L^*)^{\otimes m}\otimes O(j))=0$
for all $m$ at least equal to some $m_0=m_0(F)$. Using the injection $f$,
it follows that $H^0(X,F^*\otimes (L^*)^{\otimes m}\otimes \omega_X)=0$
for all $m\geq m_0(F)$. That is, $\Hom_X(F\otimes L^{\otimes m},\omega_X)=0$,
and hence $H^n(X,F\otimes L^{\otimes m})=0$ for all $m\geq m_0(F)$.
It follows that $L$ is $(n-1)$-T-ample, which we call simply $(n-1)$-ample
after Theorem \ref{equivred}.

Conversely, let $L$ be an $(n-1)$-ample line bundle on $X$. Then
for any line bundle $F$ on $X$, we have $H^n(X,F\otimes L^{\otimes m})=0$
for all $m$ at least equal to some $m_0=m_0(F)$. Therefore
$\Hom_X(F\otimes L^{\otimes m},\omega_X)=0$, which we write as
$H^0(X,F^*\otimes (L^*)^{\otimes m}\otimes \omega_X)=0$. Since the sheaf
$\omega_X$ is not zero, there is a nonzero map $g:O_X(-j)\arrow \omega_X$
for some $j>0$, which we fix. Here $g$ is an injection of sheaves because
the line bundle $O_X(-j)$ is torsion-free. Therefore,
$H^0(X,F^*(-j)\otimes (L^*)^{\otimes m})=0$ for all $m\geq m_0(F)$.
In particular, for every $(n-1)$-ample line bundle $L$,
$L^*$ is not big.
But the $(n-1)$-ample
cone is open in $N^1(X)$ (Theorem \ref{opencone}). So
for every $(n-1)$-ample line bundle $L$, $L^*$ is not pseudoeffective.
\qed

\section{The $q$-nef cone}
\label{q-nef}

It is an open problem to give a numerical characterization
of $q$-ampleness, analogous to the Kleiman or Nakai-Moishezon
criteria for $q=0$. We know from Theorem \ref{opencone}
that $q$-ampleness on a smooth
projective variety of characteristic zero only depends
on the numerical equivalence class of a divisor,
but that leaves the problem of describing the $q$-ample cone
by explicit inequalities. 
In this section, we disprove the most obvious attempt
at a Kleiman criterion for $q$-ampleness.

Let $X$ be a projective variety of dimension $n$ over a field.
For a natural number $q$,
define the {\it $q$-nef cone }as the set of $D\in N^1(X)$
such that for every $(q+1)$-dimensional subvariety
$Z\subset X$, $-D$ restricted to $Z$ is not big. 
The $q$-nef cone is clearly
a closed cone in $N^1(X)$, not necessarily convex. For example,
the 0-nef cone is the usual nef cone, because $-D$ is not big
on a curve $Z$ exactly when $D\cdot Z\geq 0$. Another simple case
is the $(n-1)$-nef cone, which
is the complement of the negative of the big cone in $N^1(X)$.

Another description of $q$-nef divisors comes from the theorem
of Boucksom-Demailly-Paun-Peternell.
On any complex projective variety $X$,
BDPP characterized the dual of the closed
effective cone as the closed convex cone
spanned by curves that move on $X$
\cite{BDPP}, \cite[v.~2, Theorem 11.4.19]{Lazarsfeld}.
For this statement, we say that a curve {\it moves }on a
projective variety $X$ if
it is the image under some resolution
$X'\arrow X$ of a complete intersection of ample divisors
in $X'$, $D_1\cap\cdots\cap D_{n-1}$. By BDPP's theorem,
an $\R$-divisor $D$ on $X$ is $q$-nef if and only if
for every
$q+1$-dimensional subvariety $Z$ of $X$, $D$ has nonnegative degree
on some curve that moves on $Z$.

Let $X$ be a projective variety over a field
of characteristic zero. It is straightforward to check that the $q$-ample cone 
in $N^1(X)$ is contained in the interior of the $q$-nef cone, essentially
because the restriction of a $q$-ample divisor to each $(q+1)$-dimensional
subvariety is $q$-ample. In the extreme cases $q=0$ and $q=n-1$,
the $q$-ample cone is equal to the interior of the $q$-nef cone.
But this can fail for the 1-ample cone of a smooth projective 3-fold,
as we now show. The problem remains
to give a Kleiman-type description of the $q$-ample cone.

\begin{lemma}
Let $X$ be the $\P^1$-bundle over $\P^1\times \P^1$ given
by $X=P(O\oplus O(1,-1))$, over the complex numbers.
(This is a smooth projective toric Fano 3-fold.)
Then the 1-ample cone of $X$ is strictly smaller than the interior
of the 1-nef cone of $X$.
\end{lemma}

{\bf Proof. }Let $E$ be the vector bundle $O\oplus O(1,-1)$
on $\P^1\times \P^1$. Write $X=P(E)$ for the variety of
codimension-1 subspaces of $E$, with projection
$\pi:X\arrow \P^1\times \P^1$. 
Every line bundle on $X$ has the form
$\pi^*O(a,b)\otimes O_{P(E)}(c)$ for some integers $a,b,c$.
The cohomology of a line bundle on $X$ is given by
$$H^i(X,\pi^*O(a,b)\otimes O_{P(E)}(c))=H^i(\P^1\times \P^1,
O(a,b)\otimes R\pi_*O_{P(E)}(c)).$$
If $c>0$, then $\pi_*O_{P(E)}(c)=S^c(E)=\oplus_{j=0}^c O(j,-j)$
and the higher direct images are zero. Thus, for $c>0$,
$$H^i(X,\pi^*O(a,b)\otimes O_{P(E)}(c))=H^i(\P^1\times \P^1,
\oplus_{j=0}^c O(a+j,b-j)).$$

We will compute the 1-ample cone of $X$ intersected with $c>0$,
which is enough for our purpose. If $c>0$, then $H^3(X,\pi^*O(a,b)\otimes
O_{P(E)}(c))=0$ by the previous paragraph. It follows that
$L=\pi^*O(a,b)\otimes O_{P(E)}(c)$ is 2-ample whenever $c>0$.
Next,
$$H^2(X,\pi^*O(a,b)\otimes O_{P(E)}(c))=H^2(\P^1\times \P^1,
\oplus_{j=0}^c O(a+j,b-j)).$$
Here $H^2(\P^1\times \P^1,O(p,q))$ is zero if and only if
$p\geq -1$ or $q\geq -1$.
It follows that a line bundle $L$ with $c>0$
is 1-ample if and only if $a>0$ or $b>c$ or $a+b>0$.

We now compute the 1-nef cone of $X$ (intersected with $c>0$) using
toric geometry \cite{Fulton}.
Every line bundle $L$ on $X$ can be made $T$-equivariant,
where $T\cong (G_m)^3$ acts on $X$. By definition, for $L$ to be
1-nef means that $L^*$ is not big on any surface $Y$ in $X$. Because $L$
is $T$-equivariant, if $L^*$ is big on some surface
$Y$, then $L^*$ is big on any
translate $tY$ for $t\in T$. Every action of $T$
on a projective variety has a fixed point,
and so there is a fixed point
in the closure of the $T$-orbit of $Y$ in the Hilbert scheme.
Using upper semicontinuity of $h^0$, it follows
that $L^*$ is big on some $T$-invariant subscheme $Z$ of dimension 2,
and hence on some irreducible $T$-invariant surface. Thus
$L$ is 1-nef if and only if $L^*$ is not
big on the finitely many $T$-surfaces in $X$. (This
argument shows that the $q$-nef cone of a toric variety is rational
polyhedral,
for any $q\geq 0$.)

The toric surfaces in $X$ are the two sections $S_1=P(O)$ and 
$S_2=P(O(1,-1))$ of $X\arrow \P^1\times \P^1$ and the inverse images
$Y_1,\ldots,Y_4$
of the four curves $\P^1\times 0$, $\P^1\times \infty$, $0\times \P^1$,
and $\infty\times \P^1$ in $\P^1\times \P^1$. Let $L=\pi^*O(a,b)
\otimes O_{P(E)}(c)$ be a line bundle with $c>0$. Then $L$ has positive
degree on each fiber $\cong \P^1$, and these curves cover the surfaces
$Y_1,\ldots,Y_4$. It follows that $L^*$ cannot be big on $Y_1,\ldots,Y_4$.
Thus $L$ is 1-nef if and only if $L^*$ is not big on $S_1$ and not big
on $S_2$. Here $O_{P(E)}(1)$ restricts to the trivial bundle on $S_1$
and to $O(1,-1)$ on $S_2$, so $L$ restricts to $O(a,b)$ on $S_1\cong
\P^1\times \P^1$ and to $O(a+c,b-c)$ on $S_2\cong \P^1\times \P^1$.
The big cone of $\P^1\times \P^1$ consists of the line bundles
$O(a,b)$ with $a>0$ and $b>0$. We conclude that a line bundle $L$
with $c>0$ is 1-nef if and only
if ($a\geq 0$ or $b\geq 0$) and ($a+c\geq 0$ or $b-c\geq 0$).

For example, the line bundle $\pi^*O(-2,1)\otimes O_{P(E)}(3)$ 
is in the interior
of the 1-nef cone but is not 1-ample. \qed

\section{Questions}

We raise two questions. The first is about the relation between $q$-ampleness
and K\"uronya's asymptotic cohomological functions \cite{Kuronya}.
Let $X$ be a projective
variety of dimension $n$
over a field of characteristic zero. For a line bundle $L$
on $X$, define 
$$\widehat{h}^i(L)=\limsup_{m\arrow\infty} \frac{h^i(X,L^{\otimes m})}
{m^n/n!}.$$
K\"uronya showed that $\widehat{h^i}$ extends uniquely
to a continuous function on $N^1(X)$ which is homogeneous
of degree $n$. For example, for a nef $\R$-divisor $D$, we have
$\widehat{h}^i(D)=0$ for $i>0$, and $\widehat{h}^0(D)$ is the
intersection number $D^n$.

{\bf Question. }Let $D$ be an $\R$-divisor on a projective
variety $X$ over a field of characteristic
zero. Let $q$ be a natural number.
Suppose that $\widehat{h}^i(E)$ is zero
for all $i>q$ and all $\R$-divisors $E$
in some neighborhood of $D$ in $N^1(X)$.
Is $D$ $q$-ample?

The converse is clear, since the $q$-ample cone is open in $N^1(X)$
(Theorem \ref{opencone}). The question has a positive answer 
for $q=0$, by de Fernex, K\"uronya, and Lazarsfeld \cite{DKL}.
It is also true for $q=n-1$, using that the $(n-1)$-ample cone
is the complement of the negative of the closed effective cone
(Theorem \ref{n-1}).

Another question was raised by Dawei Chen and Rob Lazarsfeld.

{\bf Question. }Let $X$ be a Fano variety, say over a field
of characteristic zero. Let $q$ be a natural number.
Is the $q$-ample cone in $N^1(X)$ the interior of a finite union
of rational polyhedral convex cones?

The answer is positive for $q=0$: the nef cone of a Fano
variety is rational polyhedral, by the Cone theorem
\cite[Theorem 3.7]{KMbook}.
It is also true for $q=n-1$, since 
the closed effective cone of a Fano variety is rational polyhedral
by Birkar, Cascini, Hacon, and M\textsuperscript{c}Kernan
\cite[Corollary 1.3.1]{BCHM}. So the first open case is
the 1-ample cone of a Fano 3-fold.


\small \sc DPMMS, Wilberforce Road,
Cambridge CB3 0WB, England

b.totaro@dpmms.cam.ac.uk
\end{document}